\documentclass[12pt,english]{article}
\usepackage{amscd,amsfonts,amsmath,amssymb,amstext,amsthm}
\usepackage{mathpazo}
\usepackage{hyperref}
\makeatletter
\renewenvironment{thebibliography}[1]
    {\noindent\\\textbf{Bibliography}
      \list{\@biblabel{\@arabic\c@enumiv}}%
           {\settowidth\labelwidth{\@biblabel{#1}}%
            \leftmargin\labelwidth
            \advance\leftmargin\labelsep
            \@openbib@code
            \usecounter{enumiv}%
            \let\p@enumiv\@empty
            \renewcommand\theenumiv{\@arabic\c@enumiv}}%
      \sloppy
      \clubpenalty4000
      \@clubpenalty \clubpenalty
      \widowpenalty4000%
      \sfcode`\.\@m}
     {\def\@noitemerr
       {\@latex@warning{Empty `thebibliography' environment}}%
      \endlist}
\makeatother
\newcommand{\publishedwork}{\\[.5cm]\hspace*{-\leftmargin}%
\textbf{Publications of Karl Stein}}
\newcommand{\dissertations}{\\[.5cm]\hspace*{-\leftmargin}%
\textbf{Dissertations guided by Karl Stein}}
\title{Karl Stein (1913-2000)}
\author{Alan Huckleberry} 
\date{}
\begin{document}
\maketitle
%

\noindent
Karl Stein was born on the first of January 1913 in Hamm
in Westfalen, grew up there, received
his Abitur in 1932 and immediately thereafter began
his studies in M\"unster. Just four years later, under
the guidance of Heinrich Behnke, he passed his Staatsexam,
received his promotion and became Behnke's assistant.

\bigskip\noindent
Throughout his life, complex analysis,
primarily in higher dimensions (``mehrere Ver\"anderliche''),
was the leitmotif of Stein's mathematics. As a fresh Ph.D. in
M\"unster in 1936, under the leadership of the master Behnke, 
he had already been exposed to the fascinating developments
in this area. The brilliant young Peter Thullen was proving 
fundamental theorems, Henri Cartan had visited M\"unster, 
and Behnke and Thullen had just written 
\emph{the} book on the subject. It must have been clear to 
Stein that this was the way to go.

\bigskip\noindent
Indeed it was! The amazing phenomenon of analytic continuation
in higher dimensions had already been exemplified more
than 20 years before in the works of Hartogs and E. E. Levi.
Thullen's recent work had gone much further. In the
opposite direction, Cartan and Thullen had proved their
characterization of domains in $\mathbb C^n$ which admit
a holomorphic function which can not be continued any further.
Behnke himself was also an active participant in mathematics
research, always bringing new ideas to M\"unster.  This was
indeed an exciting time for the young researcher, Karl Stein. 

\bigskip\noindent
Even though the pest of the Third Reich was already 
invading academia, Behnke kept things going for as 
long as possible.  But this phase of the
M\"unster school of complex analysis could not
go on forever. Although Stein was taken into the army,
during a brief stay at home he was able to prepare
and submit the paper which contained the results from
his Habilitationsarbeit which was accepted in 1940.
At a certain point he was sent to the eastern front. 
Luckily, however, the authorities were informed 
of his mathematical abilities, and he
was called back to Berlin to work until the end
of the war in some form of cryptology.  Stein told me he 
was not very good at this.

\bigskip\noindent
Almost immediately after the war, in a setting of total 
destruction, Behnke began to rebuild his group,
and very soon Stein became the mathematics guru in M\"unster. 
At the time there were only two professor positions in
pure mathematics, those of Behnke and F. K. Schmidt.  
Although it must have been very difficult, Behnke somehow found
a position for Stein which he held from 1946 und
1955.  

\bigskip\noindent
In 1955 Stein took a chair of mathematics at the 
Ludwigs-Maximilian-Universit\"at in M\"unchen where
he stayed for the remainder of his academic career.  There he
continued his mathematics and built his own group in complex analysis.
A number of his doctoral students later became professors at
universities here in Germany.
One of the most exciting periods in M\"unchen 
was certainly that in the late-1960s with the young Otto Forster,
who received his doctorate in 1961, 
leading a group of up-and-coming researchers. 

\bigskip\noindent
Not only being an outstanding researcher and teacher,
Karl Stein worked tirelessly on all sides of academia.
Among other activities he was managing editor of Manuscripta
Mathematica from 1969 until 1983, and in 1966 he was president
of the DMV.  He was awarded numerous honors, including
membership in the Bavarian and the 
Austrian Academies of Sciences, and corresponding membership 
of the G\"ottingen Academy of Sciences. In 1973 he received
an honorary doctor's degree from the faculty of mathemetics
in M\"unster, and in 1990, on the
occasion of the 100th anniversary of the founding of
the DMV, he was awarded the inaugural Cantor-Medaille.

\bigskip\noindent
Up until a few years before his death in October of 2000
Stein was still actively thinking about and even doing
mathematics. I remember his talk in Bochum in the
fall of 1992, just before his 80th birthday.
He still radiated his intense
interest in discovery and the joy of being involved
with something so beautiful. Even the youngest of students
who heard that talk were mesmerized, knowing they
had experienced the real thing! 

\bigskip\noindent
As the reader has certainly noticed we have barely touched upon the
mathematics that so fascinated Stein and his 
contributions as a researcher and teacher. Let us
devote the remainder of this article to a chronological sketch 
of some of the high points. 

\bigskip\noindent
Although Stein's thesis does not reflect his later work,
it does reflect one of the main directions of that
time, namely ``analytic continuation'', and it also
shows that even at this beginning stage he was ahead
of his time. It was already known that
a function which is holomorphic in a neighborhood of the
standard Euclidean sphere in $\mathbb C^n$, $n>1$, extends
holomorphically to the full Euclidean ball.  In his thesis
(see \cite {S1}),
under assumptions, e.g., on dimension, which we now
know to be inessential, Stein shows that such results are in fact
local in nature.  For example, a function which is holomorphic
in a neighborhood of a piece of the sphere extends to an
open set which only depends on that piece.  He even realized
that such results are possible for functions holomorphic
in neighborhoods of higher-codimensional real manifolds.
These results, which represent a change in viewpoint, 
are precursors to the highly developed modern
theory of Cauchy-Riemann manifolds.

\bigskip\noindent
One group of leading problems of that period revolved around
the question of whether or not holomorphic or meromorphic
functions could be constructed with certain prescribed properties.
The model situations were the theorem of Mittag-Leffler 
and the Weierstrass-theory of infinite product expansions
on the complex plane.  In the former case, at each
point of a divergent sequence $\{z_n\}$  
a finite negative part $P_n$ of a Laurent series is given
and one asks if there is a meromorphic function $f$ on 
the complex plane which is holomorphic everywhere except
at points of the sequence with $f-P_n$ being holomorphic near
each $z_n$.  Formulated without the details, one asks
if one can arbitrarily prescribe the principal parts
of a meromorphic function.  

\bigskip\noindent
In the original Weierstrass-theory one
prescribes an positive integer $m_n$ at each of the points
$z_n$ and asks for the existence of a holomorphic function $f$
whose zeros only occur at points of the
sequence and the orders of the zeros $f$ at these points
should be the given integers. More generally one allows
$m_n$ to be an arbitrary integer and asks for a meromorphic
function with prescribed zeros and poles. 
In this case the ``principal part'' $P_n$ is replaced
by $D_n=(z-z_n)^{m_n}$ and the requirement is that
$\frac{f}{D_n}$ is holomorphic near $z_n$. Briefly stated,
one asks if the ``divisor'' of a meromorphic function
can be arbitrarily prescribed.  

\bigskip\noindent
Due to the early work of P. Cousin (\cite{C}) one referred
to the higher-dimensional versions of these as the 
additive and multiplicative Cousin problems or simply Cousin I. 
and II.

\bigskip\noindent
As Stein was starting out, it was well-known that the 
appropriate domains  
for solving the interesting problems of the time, such as the
Cousin problems, were the
``Regularit\"atsbereiche''. Precisely speaking, they can be
defined as domains
$D$ in $\mathbb C^n$ so that given any divergent sequence 
$\{x_n\}$ in $D$ there exists a function $f$ holomorphic
on $D$ with $\mathrm {lim}\vert f(x_n)\vert =\infty$. In fact
such a domain possesses a holomorphic function which cannot
be continued across any boundary point. In other words 
$D$ is the ``region of regularity'' for that function
or its ``domain of holomorphy''. In the mid-1930s
Cartan (\cite {Ca}) and Oka (\cite {O}) had already
proved definitive results
for Cousin I for domains in $\mathbb C^n$: If $D$ is a domain of 
holomorphy, then every Cousin I problem on $D$ is solvable!

\bigskip\noindent
Immediately after his thesis Stein turned to the 
Cousin problems. Later he discovered the correct abstract
setting for solving these and many other problems, e.g., on 
complex manifolds or even complex spaces, 
but at this point his attention was focused on Cousin II
for domains in $\mathbb C^n$.

\bigskip\noindent
The situation at the time of Stein's entry into the subject
is beautifully described in (\cite {S2}).  There were already
a number of fascinating examples which showed that solving
this multiplicative problem on $D$ required more than $D$
just being a domain of holomorphy.  There was a natural way
to logarithmically change this to the additive problem,
i.e., to Cousin I, but in the process problems of well-definedness
arise.  This was not unknown in complex analysis. Monodromy, 
something in the fundamental group or first homology, was well-known,
but the obstruction to Cousin II was clearly higher order.  Nowadays we
know that this is the Chern class of the line bundle associated
to the divisor and, at least when the ambient manifold is
compact, we can regard it as the Poincar\'e dual of the
divisor itself.  But in those days these concepts
were not available.  Furthermore, even had they been on hand,
in the noncompact setting which is appropriate for Cousin II, 
relating a deRham- or Cech-class to something geometric
is not a simple matter. 

\bigskip\noindent
In the late 1930s, without modern topological
methods, but armed with strong geometric insight, this is exactly
what Stein had in mind: understanding this geometric obstruction.  
Being able to spend the year 1938
with Seifert in Heidelberg was in this regard certainly his good fortune
or maybe even fate.   In any case he returned to M\"unster
being one of the few (perhaps the only) complex analyst
who was in the position of applying ``modern'' topological
methods to problems such as Cousin II.  

\bigskip\noindent
In the work (\cite {S10}), which should be regarded as one of the most 
important in this early phase of several complex variables,
Stein completely solved Cousin II and the related
Poincar\'e problem using methods which opened doors to 
important new directions. 
The Oka principle, that a well-formulated problem in the complex 
analytic setting has a holomorphic solution on a domain of holomorphy 
if and only if it has a topological solution, could be seen
in precise form in the hands of Stein. In brief, 
modulo details which are now well-understood, here is what Stein did.

\bigskip\noindent
In its simplest form Cousin II amounts to the following:
On a domain of holomorphy $D$ we are given a 1-codimensional
subvariety $M$, i.e., a closed subset which is \emph{locally} defined
as the 0-set of a holomorphic function.  We ask for a function
which is \emph{globally} defined and holomorphic on $D$, which
vanishes exactly on $M$ and vanishes there exactly of order one.
Carefully worrying about triangulations, orientations
and all other matters that were known to be delicate in the 
infantile state of the topology of the days, he developed a 
theory which led to well-defined intersection numbers $M.K$, 
where $M$ is as above, or more generally a divisor in $D$, and 
$K$ runs through the 2-dimensional homology cycles.  
Under minor technical conditions,
even for domains finitely spread over domains of holomorphy, he
showed that a given divisor is the divisor of a meromorphic
function if and only if all of these (topologically defined!) 
intersection numbers vanish.  Not only did Stein prove this,
he could \emph{see} the topological obstruction!$\,--\,$I was
fortunate to talk with him about this on a number of occasions.  
As was mentioned above, nowadays we often only mouth something about
the Chern class, either deRahm or Cech, of the associated bundle,
and maybe we are not nearly seeing as much as Stein did in
the late 1930s!

\bigskip\noindent
Stein's, and also Behnke's, interests in Cousin type problems
were not only restricted to the higher-dimensional setting.
Although the questions they were discussing for domains in 
$\mathbb C^n$, $n\ge 2$, had long before been completely
handled for domains in the complex plane, not much was
known for general noncompact Riemann surfaces.  On the one
hand, that situation was simpler, because there were
no higher order topological obstructions.  On the other hand,
the complex analysis looked quite difficult:  Why should a 
noncompact Riemann surface possess even one nonconstant 
holomorphic function? In fact, the likes of Koebe and
Caratheodory had attempted without success to construct
such functions!

\bigskip\noindent
From their experience with higher-dimensional domains,
and knowledge of proofs of theorems of Mittag-Leffler type
for plane domains, Behnke and Stein at least knew what to 
try to do: Extend the Runge approximation theorem 
to noncompact Riemann surfaces and show that a noncompact 
Riemann surface possesses a Runge exhaustion! 
The Runge condition can be described
as follows: Let $\{U_n\}$ be an increasing sequence of open,
relatively compact subsets which exhaust the Riemann surface
$X$. Denote by $K_n$ the topological closure of $U_n$.  The exhaustion
is said to be Runge if for every $n$ every function 
holomorphic in a neighborhood of $K_n$ can be abitrarily 
well approximated in the sup-norm of $K_n$ by functions which 
are holomorphic on $U_{n+1}$. At the time it was well-known
that, e.g., for plane domains the condition that $U_n$ is
Runge in $U_{n+1}$ is equivalent to the topological condition
that the $U_n$ is relatively simply-connected in $U_{n+1}$.
In (\cite {S11}) Behnke and Stein 
succeeded in proving this
in the more general setting, thus proving
that a noncompact Riemann surface possesses a Runge 
exhaustion and as a consequence it follows that both 
Cousin I and II (\cite {S14}) have positive answers 
in that context. Due to the war-time conditions this work was published
long after its completion.

\bigskip\noindent
Up until the early 1950s Stein was still focused on the 
Cousin problems, particularly Cousin II. His last work
in this direction (\cite {S15}) may have turned out to be his 
most famous.  From this work one sees that Stein has studied
the deep and perhaps mysterious work of Oka, whom he credits 
with the theorem that on a domain of holomorphy a Cousin II 
problem is holomorphically solvable if and only if it is 
topologically solvable. 

\bigskip\noindent
As mentioned above, under a certain assumption which
would seem only to be technical, Stein had made this precise
in terms of his intersection numbers.  This assumption is
that the first homology group of the domain should have
a basis.  Here Stein observes that (believe it or not!)
this is really an assumption, and in order to do away with
it he must refine his topological condition.  Underway
he even proves several new results for countable 
Abelian groups!

\bigskip\noindent
Of course (\cite {S15}) is a basic work, but the reason that
it may be one of Stein's most famous is that, without pursuing
matters much further, he noted that most results of the
type he had been considering are true for, in Stein's words
and notation, 
domains $\mathfrak {G}$ in complex manifolds $\mathfrak {M}^{2n}$ which 
satisfy the following three axioms:
\begin {enumerate}
\item 
{\bf (Holomorphic convexity)} For every compact subset 
$\mathfrak {G}_0$ of $\mathfrak {G}$ there is a compact 
subset $\mathfrak {G}_1$ which contains it so that 
for every point $P$ in $\mathfrak {G}$ which is not in 
$\mathfrak {G}_1$ there is a holomorphic function $f_P$ on 
$\mathfrak {G}$ with
$$
\vert f_P(P)\vert \ >\ \mathrm{Max}\vert f_P(K_0)\vert \,.
$$
\item
{\bf (Point separation)} For any two different points 
$P_1$ and $P_2$ in $\mathfrak {G}$ there is a function $f_{P_1,P_2}$
which is holomorphic on $\mathfrak {G}$ and which takes on
different values at $P_1$ and $P_2$.
\item
{\bf (Coordinates)} For every $Q$ in $\mathfrak {G}$ there is a system
of $n$ holomorphic functions on $\mathfrak {G}$ whose functional 
determinant at $Q$ is nonzero.
\end {enumerate}
The Cartan-Serre theory, in particular the vanishing
theorems for cohomology defined by coherent sheaves 
on spaces which satisfy these axioms, was announced
by Cartan at the famous \emph{Colloque sur les fonctions
de plusieurs variables} in Brussels in 1953.  There
he baptized these spaces \emph{Variet\'e de Stein}, a
notation that is still used today. During my
very first seminar talk where Stein was present, \emph{his} manifolds
arose and, noticing my nervousness, without prompting, he
said, ``I like to call them holomorphically complete''.

\bigskip\noindent
Returning to M\"unster after participating in the Brussels
Colloquium where he announced his own fundamental work
on analytic decompositions, Stein lamented, 
``Die Franzosen haben Panzer, wir nur Pfeile und Bogen''
\footnote {Oral communcation from R. Remmert.  See (\cite {R})
for other recollections of the spirit of those
times.} To a certain extent this analogy might fit, but 
in appearance only .  Looking back one sees that
these ``Bows and Arrows'' were really quite sophisticated and
that the accomplishments of the M\"unsteraner 
were truly extraordinary!

\bigskip\noindent
The most well-known names associated with the
early days of the postwar M\"unster school of
Heinrich Behnke are Hirzebruch, Grauert, Remmert
and Stein.  
Hirzebruch, who was one of the first doctoral
students after the war, went on to prove numerous important 
results in complex geometry, primarily for compact manifolds. 
Certain of his fundamental works utilize topological methods 
which go well beyond those employed by Stein, but which are of
a similar basic spirit in that invariants such as characteristic 
classes or intersection numbers are fundamental topological 
obstructions to solving problems of analytic or algebraic geometric 
interest.  In the early days he and Stein often commuted together
from Hamm (Hirzebruch also grew up there), sometimes having to
ride on the outide running board of the train, but nevertheless 
discussing mathematics.  I can imagine that Stein's animated
expositions about his intersection numbers, or whether or not 
the first Betti group has a basis, made a lasting impression on 
the young Hirzebruch!

\bigskip\noindent
Certain of Grauert's early works, e.g., his Oka principle,
can be regarded as taking Stein's prewar mathematics 
to another universe (see, e.g., our article, \emph{Hans Grauert:
Mathematician Pur, Mitteilung of the DMV, 2008}, for a brief
summary of Grauert's work).  Later on (Stein had been retired for
a number of years) they had close common interests in 
understanding the conditions under which the quotient of 
a complex space by an analytic or meromorphic equivalence 
relation is again a complex space. I recall several
\emph{very} animated discussions in Oberwolfach!

\bigskip\noindent
In any account of Stein's mathematics after his period
of intense interest in the Cousin problems, in particular
in the topological obstructions, his work 
with Reinhold Remmert must have center stage. This turned the
page to a completely new direction!

\bigskip\noindent
Very early in Remmert's studies, Behnke sent him to 
Stein, who at the time had an idea that analytic continuation 
was something that applied not only to functions.  Maybe 
Thullen's result in the 1-codimensional case could be
proved for general analytic sets!  Stein had in mind that the
appropriate \emph{elimination theory} could be
found in Osgood's book and Remmert should check this. 
What a daunting task for someone just starting out!  
As it turned out, nothing of this sort could be found in 
Osgood, and work could be started
toward what would be the Remmert-Stein extension 
theorem (\cite {S18}).

\bigskip\noindent
Here is a statement of the simplest version of that result: 
Let $E$ be an analytic set in a domain $D$ in $\mathbb C^n$, i.e., 
a closed subset which is locally defined as the common 0-set of finitely
many holomorphic functions, and suppose that $A$ is an analytic set
in the complement $D\setminus E$ which is everywhere of
larger dimension than $E$.  Then the topological
closure $\bar A$ of $A$ in $D$ is an analytic subset of $D$
and what one adds to $A$ to obtain this closure is just
the lower-dimensional analytic subset $\bar A\cap E$.

\bigskip\noindent
To the ear of the nonspecialist the above may sound overly 
complicated. However, considering the following example,
which was a starting point for the Remmert-Stein
discussions, should allay any doubts about its importance. 
Let $D$ be $\mathbb C^n$ itself and $E$ just be the origin.
Assuming that $A$ is everywhere at least 1-dimensional,
in this case the theorem just says that $\bar A=A\cup \{0\}$
is an analytic subset of $\mathbb C^n$ and, using
results that were already known at the time, $\bar A$
is the common 0-set of finitely many holomorphic functions
which are \emph{globally} defined on $\mathbb C^n$, i.e.,
convergent power series.

\bigskip\noindent
Preimages $A=\pi ^{-1}(V)$ via the standard projection 
$\pi : \mathbb C^n\setminus \{0\}\to \mathbb P_{n-1}(\mathbb C)$
of analytic sets $V$ in projective space are examples
of analytic sets where the Remmert-Stein theorem can be applied.
In this case $A$ is invariant by the $\mathbb C^*$-action defined
by scalar multiplication.  Thus, writing the defining
power series $A$ as sums of homogeneous terms, one shows that
$A$ is also the common 0-set of finitely many of these 
\emph{homogeneous polynomials}. Consequently the original variety $V$ 
is the common 0-set of the same polynomials and is therefore 
an algebraic variety.

\bigskip\noindent
The above proof of Chow's theorem was given \emph{ahead of
time} by Cartan in his lecture at the International Congress 
of Mathematicians in Boston in 1950!  This result is 
a first example of a general principle which states
that in many algebraic geometric settings there is
no difference between algebraic and analytic phenomena.
The Remmert-Stein theorem is certainly one of the guiding
forces behind this principle!

\bigskip\noindent
The theme of holomorphic and meromorphic maps was
one of Stein's favorites and throughout this area
the Remmert-Stein theorem plays a key role. The
idea, e.g., for analyzing a holomorphic map $F:X\to Y$, is to throw out
the analytic subsets (images and preimages) where
$F$ degenerates, prove a good result for the restricted
map, and then obtain the desired result by Remmert-Stein
continuation.  In several complex variables,
meromorphic maps have indeterminacies and thus 
it is necessary to define such via their
graphs. In any theory for these \emph{set valued maps}
the Remmert-Stein result is used at many steps along the way.
Remmert developed this theory for (generically single-valued)
meromorphic maps, and Stein later generalized this
to \emph{correspondences} which are
not necessarily generically single-valued (see, e.g.
\cite{S34,S35}).

\bigskip\noindent
Remmert's mapping theorem, \emph{Images of analytic sets 
under proper holomorphic maps are analytic sets}, is very 
much in the spirit of the times. Of course this result
is extremely useful. However, it is perhaps just as important
that it calls our attention to the concept ``proper'',
i.e., inverse images of compact sets are compact.  Its role 
had already been emphasized by Henri Cartan in 1935 in the
context of actions automorphism groups on bounded domains
and some basic results were proved in Bourbaki, but
the proper mapping theorem and Stein's fundamental paper
on analytic decompositions (\cite {S23}) cemented the position 
of properness in complex analysis.

\bigskip\noindent
Stein's paper contains a wealth of interesting and useful results,
some even at the general topological level (see for example
Satz 9), but due to lack of space we will only extract
the most well-known one.  For this it should be recalled
that, in M\"unster, complex spaces were defined as topological
spaces which could be locally realized as finite ramified
covers (with obvious topological assumptions) over domains in 
$\mathbb C^n$. Stein had in fact shown that unramified
(even infinite) covers of holomorphically complete spaces 
are holomorphically complete (\cite {S24}), but he had really 
focused his interests on situations where some sort of
properness is available.   

\bigskip\noindent
Let us state an example of a result which is an important special
case of those in (\cite {S23}).  Suppose $F:X\to Y$ is a proper
holomorphic mapping of complex spaces. The domain space $X$
is assumed to be normal$\,--\,$for our purposes here it is enough
to consider the smooth case. In order to analyze
$F$, first apply Remmert's theorem so that it may be assumed
that it is surjective. Then define an equivalence relation 
$\sim $ on $X$ with two points being equivalent whenever 
they are in the same connected component of an $F$-fiber.  
The decomposition of $X$ into equivalence classes is a 
special case of what Stein called an ``analytic decomposition''.
In this case at hand, he shows that $X/\sim =:X^*$ carries a 
unique structure of a normal complex space such that 
the quotient map $\Phi:X\to X^*$ is holomorphic and every
other holomorphic map which is constant on the equivalence
classes of $\sim $ factors through it.  In particular,
this induces a holomorphic map $f:X^*\to Y$ which is a finite
ramified cover!  The factorization $F=f\circ \Phi $ is
what is now called the \emph{Stein factorization} of $F$. 

\bigskip\noindent
A number of Stein's last published works are devoted to understanding
more general situations where it is possible to construct a universal
quotient of the above type. The works (\cite {S29,S30}) are
typical of this. One exception is (\cite {S27}).
In this jewel, given two (concrete) domains in $\mathbb C^n$, 
Remmert and Stein study the possibilities for
proper holomophic maps between them.  For two polyhedral domains 
$A$ and $A^*$ with sufficient structure coming from the
affine structure of $\mathbb C^n$, they show that proper
holomorphic maps which respect this structure are in fact
affine.  In particular, for domains in $\mathbb C^2$ this
leads to strong nonexistence (rigidiy) results, e.g.,
that certain very simple explicitly given domains have
only the identity as proper holomorphic self-maps.
Their methods even shed new light on situations which were 
classically ``understood''.
For example, Poincar\'e showed that the Euclidean ball
$B_2:=\{(z,w)\in \mathbb C^2; 
\vert z\vert ^2\ + \ \vert w\vert ^2<1\}$ and the polydisk 
$\Delta_2:=\{(z,w)\in \mathbb C^2;
\vert z\vert <1\ \text{and} \ \vert w\vert < 1\}$ are not
equivalent by a biholomorphic map, because their automorphism
groups don't have the same dimensions.  
Remmert and Stein show that, just as the beginner would like to believe, 
the reason for the inequivalence of these domains is that the 
boundary of $B_2$ is round and most of the boundary of 
$\Delta $ is flat!

\bigskip\noindent
We have now come to the end of our tour of what we find to be 
the highest points of Karl Stein's mathematical works and would
like to close this note by expressing our greatest
respect and admiration, not only for the science of the man,
but equally for the man behind the science!

\bigskip\noindent
Fakult\"at und Instit\"ut f\"ur Mathematik\\
Ruhr-Universit\"at Bochum\\
44780 Bochum\\
ahuck@cplx.rub.de
\end {document}